\documentclass[12pt]{article}
\usepackage{amssymb} \usepackage{amsfonts}
\usepackage{amsmath}

\newcommand{\R} {{\mathbb R}}
\newcommand{\Q} {{\mathbb Q}}
\newcommand{\Z} {{\mathbb Z}}

\newcommand{\N} {{\mathbb N}}

\newcommand{\GL}{{\mathrm{GL}}}

\author{Rafa\l {} Lutowski\thanks{Supported by University of Gda\'nsk grant number BW - 5107-5-0345-0}{ }
 Institute of Mathematics\\
University of Gda\'nsk\\
ul. Wita Stwosza 57\\
80-952 Gda\'nsk, Poland\\
E-mail: \texttt{rlutowsk@mat.ug.edu.pl}}

\title{Seven dimensional flat manifolds with cyclic holonomy}

\begin{document}

\maketitle

\begin{abstract}
We classify (up to affine equivalence) all $7$-dimensional flat manifolds
with cyclic holonomy group.
\end{abstract}

\section{Introduction}

Let $M^n$ be a closed Riemannian manifold of dimension $n.$
We shall call $M^n$ flat if, at any point, sectional curvature is equal to zero. Equivalently,$M^n$ is isometric
to the orbit space $\R^n/\Gamma,$ where $\Gamma$ is a discrete, torsion-free and co-compact subgroup of $O(n)\ltimes\R^n$ = Isom($\R^n$). 
From the first Bieberbach theorem (see \cite{charlap}, \cite{S}, \cite{wolf}) $\Gamma$ defines a~short exact sequence of groups
\begin{equation}\label{crysb}
0\rightarrow \Z^n\rightarrow\Gamma\stackrel{p}\rightarrow G\rightarrow 0,
\end{equation}
where $G$ is a finite group. 
$\Gamma$ is called a Bieberbach group and $G$ its holonomy group.
Moreover, from second and third Bieberbach theorems (see \cite{charlap}, \cite{S}, \cite{wolf})
there are only finite number of the isomorphism classes of Bieberbach groups
of dimension $n$ and two Bieberbach groups are isomorphic if and only if they are
conjugate in the group $\GL(n,\R)\ltimes\R^n.$

With support of a computer system CARAT (\cite{carat}) it is possible
to give a~complete list of all isomorphism classes of Bieberbach groups up to
dimension $6.$ Moreover for a finite group $G$ and a number $n,$ CARAT gives possibilty for a classification (up to isomorphism)
of all Bieberbach groups of a dimension $n$ with a holonomy group $G.$
In this article the CARAT system is used to
calculate a list of all isomorphism classes of $7$-dimesnional Bieberbach groups
with cyclic holonomy group. The final list of 316 groups is presented on the www page (see \cite{Lu10}),
where the method of exposition is borrowed from \cite{carat}. 
Our main motivation was a paper \cite{S1} about $\eta$-invariants of flat manifolds,
where our results are applied.
 
A holonomy representation $\phi:G\to \GL(n,\Z)$ of the Bieberbach group $\Gamma$
(cf. \eqref{crysb}) is defined by the formula:
\begin{equation}\label{holonomyrep}
\forall g\in G,\phi(g)(e_i) = \tilde{g} e_i \tilde{g}^{-1},
\end{equation}
where $e_i\in\Gamma$ are generators of $\Z^n$ for $i=1,2,...,n,$ and $\tilde{g}$
is an element of $\Gamma$, such that $p(\tilde{g})=g.$

\section{$\Q$-classes of holonomy representation}

By \cite{Hi85}, the possible orders of cyclic groups, 
that can be realized as holonomy groups of crystallographic groups 
in dimension 7 are: 1, 2, 3, 4, 5, 6, 7, 8, 9, 10, 12, 14, 15, 18, 20, 24, 30. 
By \cite[Lemma 2.1]{Hi85}, the degree of a matrix with $n$-th primitive root of 1 is not less than 
$\varphi(n),$ where $\varphi$ is the Euler's function. 
Since for $n=15,20,24,30\hskip 2mm \varphi(n) > 7$, then any matrix of order $n$ and degree $7$ 
must be taken as a~direct sum of matrices, which orders are proper divisors of $n$. 

Let $n \in \N, n > 1$ and
$$
\Phi_{n}(x) = x^{\varphi(n)} + a_{\varphi(n)-1}x^{\varphi(n)-1} + \ldots + a_1 x + a_0
$$
be the cyclotomic polynomial of order $n$ (see \cite[page 137]{CuRe62}).
Since the characteristic polynomial of the matrix
\[
A_n = 
\begin{bmatrix}
0 & 0 & \ldots & 0 & -a_0\\
1 & 0 & \ldots & 0 & -a_1\\
0 & 1 & \ldots & 0 & -a_2\\
\vdots & \vdots & \ddots & \vdots & \vdots \\
0 & 0 & \ldots & 1 & -a_{\varphi(n)-1}
\end{bmatrix}
\]
is equal to $\pm\Phi_n(x)$, then eigenvalues of $A_n$ are primitive $n$-th roots of the unity. 

Let $G$ be a~cyclic group of order $n$, generated by an element $g$. 
Then for each $d \mid n$, $\rho_d \colon G \to \GL(\varphi(d),\Z)$, given by
\[
\rho_n(g) = A_d
\]
is an integral representation of $G$, which is irreducible over $\Q$. 
Moreover, by \cite[Corollary 39.5]{CuRe62}, these are all, up to equivalence, rational irreducible representation of $G$.

From the above remarks, rational representations of a cyclic group of order $n$ in dimension 7 are of the form
\begin{equation}
\label{eq:repq}
\rho = \mathop{\bigoplus_{d\mid n}}_{d \leq 18} a_d\rho_d,
\end{equation}
where $a_i \in \N$ and
\[
\mathop{\sum_{d\mid n}}_{d \leq 18} a_d \varphi(d) = 7,
\]
and $\rho$ is faithful, if
$$
\mathrm{LCM}\{d\mid n ;  a_d \neq 0\} = n.
$$

In the Table \ref{tab:cycpols} we give a~list of cyclotomic polynomials 
for given $n$ and some remarks about the matrices $A_n$. 
The relation $\sim$ means "the same conjugacy class in $\GL(\varphi(n),\Q)$".

\begin{table}
\begin{center}
\begin{tabular}{c|c|r|l}
 $n$ & $\varphi(n)$ & $\Phi_n(x)$ & Remarks\\ \hline
 $2$ &        $1$ &                     $x-1$\\
 $3$ &        $2$ &                 $x^2+x+1$\\
 $4$ &        $2$ &                   $x^2+1$\\
 $5$ &        $4$ &         $x^4+x^3+x^2+x+1$\\
 $6$ &        $2$ &                 $x^2-x+1$ & $A_6 \sim -A_3$\\
 $7$ &        $6$ & $x^6+x^5+x^4+x^3+x^2+x+1$\\
 $8$ &        $4$ &                   $x^4+1$\\
 $9$ &        $6$ &               $x^6+x^3+1$\\
$10$ &        $4$ &         $x^4-x^3+x^2-x+1$ & $A_{10} \sim -A_5$\\
$12$ &        $4$ &               $x^4-x^2+1$\\
$14$ &        $6$ & $x^6-x^5+x^4-x^3+x^2-x+1$ & $A_{14} \sim -A_7$\\
$18$ &        $6$ &               $x^6-x^3+1$ & $A_{18} \sim -A_9$\\
\end{tabular}
\end{center}
\caption{Cyclotomic polynomials for given numbers}
\label{tab:cycpols}
\end{table}  

\section{Determination of Bieberbach groups}

Let $G$ be a~cyclic group of order $n$. 
From \eqref{eq:repq} we know, how to determine all equivalence classes of seven dimensional rational representation of $G.$ 
We want to classify  
all Bieberbach groups with a holonomy group $G.$ 
There are three steps:
\begin{enumerate}
\item Determine, up to equivalence, all faithful representations $\rho \colon G \to \GL(7, \Q);$
\item Determine all integral representations (up to equivalence) of $G$ equivalent over $\Q$ to $\rho;$
\item For each representation $\tau \colon G \to \GL(7,\Z)$ from the previous point, 
determine all Bieberbach groups (up to isomorphism) with holonomy representation $\tau$.  
\end{enumerate} 
To determine $\Z$-classes of faithful representations of cyclic 
group of prime order, we use \cite[Theorem 74.3]{CuRe62}.

As mentioned before, the complete list of seven dimensional Bieberbach groups with a cyclic holonomy group is given in \cite{Lu10}. 
Let us give a short dictionary of tables. If  a~Bieberbach group $\Gamma$ has a~name of the form
\[
\mathtt{n/n.a_1xf_1\_b_1\text{+}\ldots\text{+}a_lxf_l\_b_l.p.q.r},
\]
then $n$ is the order of the holonomy group $G,$ the $\Q$-class of a holonomy 
representation of $G$ is given by the representation 
\[
a_1\rho_{b_1}\oplus\ldots\oplus a_l\rho_{b_l}
\]
(cf. \eqref{eq:repq}). Moreover  $f_i = \varphi(b_i)$, for $i=1,\ldots,l;$ $p.q$ is a~symbol of the 
$\Z$-class of a holonomy representation and $r$ is a number of the group $\Gamma$. The numbers $p,q,r$ are assigned by CARAT.


\end{document}